\documentclass[11pt, a4paper]{article}
\usepackage{amsmath,amsfonts,amssymb,amsthm,epsfig,epstopdf}
\usepackage{geometry}
\usepackage{graphicx,graphics}
\usepackage{txfonts}
\usepackage{setspace}
\usepackage{parskip}
\usepackage{color}
\theoremstyle{definition}

\theoremstyle{plain}

\newcommand\ci{\perp\!\!\!\perp}
\usepackage{tikz}
\usepackage{pgfplots}
\pgfplotsset{compat=newest}
\usepackage{verbatim}
\usepackage{fancyvrb}
\usepackage[bookmarks=false]{hyperref}
\usetikzlibrary{arrows}
\usetikzlibrary{positioning}
\usetikzlibrary{shapes.geometric,positioning}
\usepackage{graphics, bm}
\usetikzlibrary{matrix}
\usepackage{caption}
\usepackage{subcaption}
\usepackage{amsmath,amsfonts,amssymb,amsthm,epsfig,epstopdf}
\usepackage{txfonts}
\usepackage{subcaption}
\usepackage[framemethod=TikZ]{mdframed}
\usepackage{cleveref}
\usepackage{authblk}
\usepackage{appendix}
\BeforeBeginEnvironment{appendices}{\clearpage}
\usepackage{url}
\usetikzlibrary{automata,positioning,arrows}
\usepackage{tabularx,ragged2e,booktabs,caption}
\newcolumntype{C}[1]{>{\Centering}m{#1}}

\newcommand{\rpm}{\sbox0{$1$}\sbox2{$\scriptstyle\pm$}
  \raise\dimexpr(\ht0-\ht2)/2\relax\box2 }
\mdfdefinestyle{MyFrame}{%
    linecolor=black,
    outerlinewidth=1pt,
    innertopmargin=\baselineskip,
    innerbottommargin=\baselineskip,
    innerrightmargin=20pt,
    innerleftmargin=20pt,
    backgroundcolor=white}
\usepackage[numbers, round]{natbib}
\title{Longitudinal Mediation Analysis Using Natural Effect Models}
\author[1]{Murthy N. Mittinty \thanks{murthy.mittinty@adelaide.edu.au}}
\author[2,3]{Stijn Vansteelandt \thanks{stijn.vansteelandt@ugent.be}}
\affil[1]{\footnotesize School of Public Health, The University of Adelaide, South Australia, 5000 Australia. Phone:+61 8 8303961}
\affil[2]{\footnotesize Department of Applied Mathematics, Computer Science and Statistics, Krijslaan 281, Gent University, Gent , Belgium.}
\affil[3]{\footnotesize Department of Medical Statistics, London School of Hygiene and Tropical Medicine, London, UK}
\date{}
\begin{document}
\maketitle
\begin{abstract}
\doublespacing
Mediation analysis is concerned with the decomposition of the total effect of an exposure on an outcome into the indirect effect through a given mediator, and the remaining direct effect. This is ideally done using longitudinal measurements of the mediator, as these capture the mediator process more finely. However, longitudinal measurements pose challenges for mediation analysis. This is because the mediators and outcomes measured at a given time-point can act as confounders for the association between mediators and outcomes at a later time-point; these confounders are themselves affected by the prior exposure and outcome. Such post-treatment confounding cannot be dealt with using standard methods (e.g. generalised estimating equations). Analysis is further complicated by the need for so-called cross-world counterfactuals to decompose the total effect. This article addresses these challenges. In particular, we introduce so-called natural effect models, which parameterise the direct and indirect effect of a baseline exposure w.r.t. a longitudinal mediator and outcome. These can be viewed as a generalisation of marginal structural models to enable effect decomposition. We introduce inverse probability weighting techniques for fitting these models, adjusting for (measured) time-varying confounding of the mediator-outcome association. Application of this methodology uses data from the Millennium Cohort Study, UK.
\end{abstract}

\section{Introduction}

Mediation analyses are ideally based on longitudinal measurements of the mediator. These represent the entire mediator process better than a single assessment of the mediator. At least in principle, they should thus better allow to capture the extent to which exposure influences outcome by affecting that mediator process, which may in turn influence outcome. In spite of this, developments on counterfactual-based mediation analyses have rarely considered longitudinal mediators. This is largely related to the difficulties of identifying (natural) direct and indirect effect in the presence of confounders of the mediator-outcome association, which are themselves affected by the exposure \cite{asp05,VVR}. Such confounders are essentially guaranteed to exist in longitudinal mediation analyses because the association between mediator and outcome at a given point in time is typically confounded by the history of mediators and outcomes (amongst other things), which are often affected by the exposure.

In view of this, early work on longitudinal counterfactual-based mediation analysis \cite{bi16} considered settings where such time-varying confounders could be assumed absent. More recent work \cite{vtlm17,zv18} has shifted the focus from natural direct and indirect effect (or more specifically, path-specific effects) to so-called interventional direct and indirect effect, which can be identified under weaker conditions \cite{VVR,joh17,VD17}, but do not always enable decomposition of the total effect into direct and indirect pathways. By relying on a general identification theory for path-specific effects by Shpitser \cite{sh13}, Vansteelandt et al. \cite{SV17} showed how to decompose the effect of a point exposure on an arbitrary outcome into the path-specific effect via the mediator (i.e., the combined effect due to the exposure directly influencing one of the mediators, which then in turn influences the outcome) versus the remaining direct effect. Assuming that the data-generating mechanism can be represented by a nonparametric structural equation model with independent errors, and assuming that there is no unmeasured confounding of the exposure-mediator, exposure-outcome and mediator-outcome associations, they used strategies analogous to g-computation to estimate these effects from observational data. Like other g-computation estimators, their results can be very sensitive to model misspecification, do not readily allow for testing for the presence of direct and indirect effect, and can be difficult to report when the exposure takes on more than two levels or the outcome is repeatedly measured over time (for then one may need to report a direct and indirect effect for each exposure level and w.r.t. each outcome separately).

To accommodate this, this paper extends the class of natural effect models previously introduced by Lange et al.\cite{lg12} and Vansteelandt et al.\cite{sv12} to longitudinal mediation analysis. Natural effect models generalise marginal structural models to enable effect decomposition. Building on \cite{zv18} we propose inverse probability weighting strategies for fitting these models, which are much less demanding in terms of modelling assumptions than g-computation estimators, are relatively easy to perform using standard software, and avoid extrapolation when subjects with different levels of exposure or mediator are very different in their observed background characteristics. We illustrate the proposal with an analysis of the effect of socio-economic status on child mental health using data from the Millennium Cohort Study (UK) and the extent to which this effect is mediated by maternal psychological distress.

\section{Effect decomposition into direct and indirect effect}
\subsection{Decomposition in single mediator and single exposure setting}
\textit{Notation, definition and identification}. In the counterfactual framework, causal effects are defined by contrasting counterfactual outcomes under different exposure settings. For example, the total causal effect of a binary exposure ($A=1$ for exposed, $A=0$ for unexposed) on an outcome, $Y$, is obtained by comparing $Y_1$ and $Y_0$, with $Y_a$ the counterfactual outcome that would have been observed if $A$ were set, possibly contrary to the fact, to $a$. The population average effect then can be quantified in terms of a mean difference, $E\{Y_1-Y_0\}$, or when the outcome is binary, alternatively as the relative risk $P\{Y_1=1\}/P\{Y_0=1\}$.
Following the causal inference literature \cite{tv15,Pearl,vv09}  we will further describe direct and indirect effect via a given mediator $M$ in terms of the nested counterfactuals $Y_{a,M_{a^*}}$, the counterfactual outcome that would have been observed if $A$ was set to $a$ and $M$ was set to the value it would have taken if $A$ was set to $a^{*}$. In particular we will compare $Y_{a,M_{a^*}}$ with $Y_{a^*,M_{a^*}}$ to express the direct effect of changing the exposure $a$ to $a^*$. Such comparison can for instance, be made in terms of an average difference within levels of baseline covariates ($L_0$), $E\{Y_{a,M_{a^*}}-Y_{a^*,M_{a^*}}|L_0\}$, or marginally, $E\{Y_{a,M_{a^*}}-Y_{a^*,M_{a^*}}\}$; as a risk ratio, $P\{Y_{a,M_{a^*}}=1\}/P\{Y_{a^*,M_{a^*}}=1\}$, and so on. Likewise we will compare $Y_{a^*,M_{a}}$ with $Y_{a^*,M_{a^*}}$, to obtain a measure of the indirect effect via the mediator. For example, on the additive scale, the total causal effect decomposes into the sum of the so-called natural direct effect and indirect effect
\[
E\{Y_1-Y_0\}=E\{Y_{1,M_0}-Y_{0,M_0}\}+E\{Y_{1,M_1}-Y_{1,M_0}\}
\]
given the composition assumption that $Y_{a,M_a}=Y_a$. The word ``natural" refers to the fact that we have let the mediator take the value it would take naturally when the exposure is set to $a$.

\begin{figure}
\centering
\begin{subfigure}[b]{0.475\textwidth}
\captionsetup{justification=centering}
\begin{tikzpicture}[->,>=stealth',auto,node distance=3cm,
  thick,main node/.style={circle,draw,font=\sffamily\Large\bfseries}]

  \node[main node] (1) {$A$};
  \node[main node] (2) [right of =1]{$M$};
  \node[main node] (3) [above of =2]{$L_0$};
  \node[main node] (4) [right of =2] {$Y$};

  \path[every node/.style={font=\sffamily\small}]
    (3) edge node [below][] {} (1)
    (3) edge node [below][] {} (2)
    (3) edge node [below][] {} (4)
    (2) edge node [right] [midway, above]{$II$} (4)
    (1) edge[bend right] node [right] [midway, above] {$I$} (4)
    (1) edge node [right][midway,above] [midway,above] {$II$} (2);
\end{tikzpicture}
%\caption{\small{\label{first}{Mediation Analysis with Single Mediator}}}
\end{subfigure}
\begin{subfigure}[b]{0.475\textwidth}
\captionsetup{justification=centering}
\begin{tikzpicture}[->,>=stealth',auto,node distance=3cm,
  thick,main node/.style={circle,draw,font=\sffamily\Large\bfseries}]

  \node[main node] (1) {$A$};
  \node[main node] (2) [right of =1]{$M$};
  \node[main node] (3) [above of =2]{$L_1$};
  \node[main node] (4) [right of =2] {$Y$};

  \path[every node/.style={font=\sffamily\small}]
    (1) edge node [above][sloped,anchor=center, above] {$I$} (3)
    (3) edge node [below][sloped, anchor=center, above] {$I$} (2)
    (3) edge node [below][sloped,anchor=center, above] {$I$} (4)
    (2) edge node [right] [midway, above]{$II$} (4)
    (1) edge[bend right] node [right] [midway, above] {$I$} (4)
    (1) edge node [right][midway,above] [midway,above] {$II$} (2);
\end{tikzpicture}
%\caption{\small{Mediation with exposure-induced-mediator-outcome confounder}}
%\label{f1b}
\end{subfigure}
\caption{{\small Causal Directed Acyclic Graph with exposure A, mediator M, outcome Y and pre-treatment confounder $L_0$ (left) and post treatment confounder $L_1$ (right)}}
\label{f1}
\end{figure}
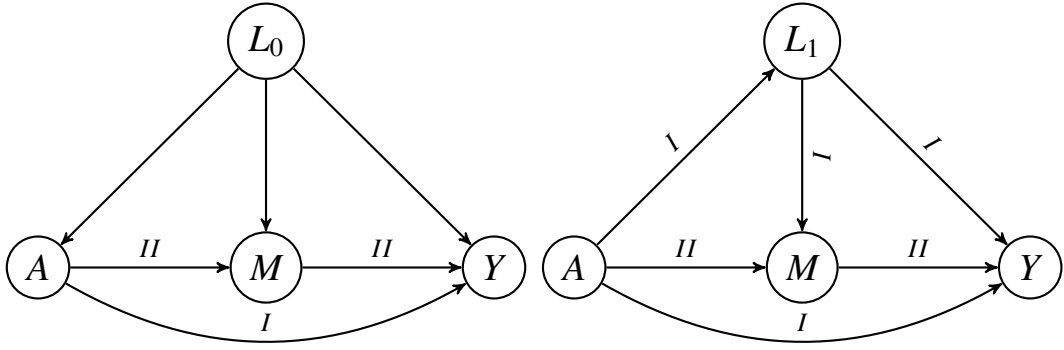

In this one time point setting, nonparametric identification of natural direct and indirect effect is possible by making a set of sufficient conditions, which state that for any value of $a,a^*$ and $m$
\begin{equation}\label{eq1}
Y_{a,m}\ci A|L_0
\end{equation}
\begin{equation}\label{eq2}
Y_{a,m}\ci M|A=a,L_0
\end{equation}
\begin{equation}\label{eq3}
M_{a}\ci A|L_0
\end{equation}
\begin{equation}\label{eq4}
Y_{a,m}\ci M_{a^*}|L_0
\end{equation}
where $A\ci B|C$ must be read as $A$ and $B$ are independent conditional on $C$. Here, $Y_{am}$ denotes the counterfactual outcome that would have been observed if $A$ were set to $a$ and $M$ to $m$.
Conditions 1-4 require just the baseline confounders $L_0$ to deconfound (a) the effect of exposure $A$ on outcome $Y$, (b) the effect of mediator $M$ on outcome $Y$ conditional on exposure; and (c) the effect of the exposure on the mediator.  Assumption (4) is stronger because it involves the dependence between counterfactuals at different exposure levels. Assumptions (1)-(4) cannot simultaneously hold when there are confounders of the mediator-outcome association that are affected by the exposure.\\
\textit{Natural effect models}. Under conditions 1-4, the natural direct and indirect effect can be parameterised via so-called natural effect models  \cite{lg12,sv12}. These express the mean of nested counterfactual outcomes, thereby naturally extending marginal structural models to allow for effect decomposition. For example, suppose that the mean of $Y_{a,M_{a^*}}$ obeys
\begin{equation}\label{eq5}
E(Y_{a,M_{a^*}})=\beta_0+\beta_{1} a+\beta_{2} a^*+\beta_{3} a.a^*
\end{equation}
for all $a, a^*$. Then, $\beta_1$ captures the direct effect
\[
E\{Y_{1,M_0}-Y_{0,M_0}\}=\beta_1
\]
and $\beta_2+\beta_3$ captures the indirect effect
\[
E\{Y_{1,M_1}-Y_{1,M_0}\}=\beta_2+\beta_3
\]
The effects $\beta_1$ and $\beta_2+\beta_3$ correspond to the $A\rightarrow Y$ and the $A\rightarrow M\rightarrow Y$ paths, labelled as $I$ and $II$ in the directed acyclic graph (DAG) (see Figure \ref{f1}, left panel).
This decomposition of the total causal effect $(\beta_1+\beta_2+\beta_3)$ is not unique \cite{Greenland}. In particular, the direct effect can alternatively be defined as
\[
E\{Y_{1,M_{1}}-Y_{0,M_{1}}\}=\beta_1+\beta_3
\]
and the indirect effect as
\[
E\{Y_{0,M_1}-Y_{0,M_0}\}=\beta_2
\]
It is thus seen that when $\beta_3 \ne 0$, the direct effect may depend on the natural level at which the mediator is controlled.
Finally, note that when $a=a^*$, then model \ref{eq5} reduces to the marginal structural model
\[E(Y_{a,M_a})=E(Y_{a})=\beta_0+(\beta_{1} +\beta_{2} +\beta_{3}a) a.\]
Natural effect models thus generalise marginal structural models to enable effect decomposition.

Model \ref{eq5} is a special case of the wider class of generalised linear natural effect models \cite{lg12,sv12} given by
\begin{equation*}
g\{E(Y_{a,M_{a^*}}|L_{0})\}=\beta^{\mathsf{T}}W(a,a^*,L_{0}^*)
\end{equation*}
where $g$ is a known link function (e.g. the identity, log or logit link) and $W(a,a^*,L_0^*)$ is a known vector with components that may depend on $a$ and $a^*$ and (possibly) a set of baseline covariates $L_{0}^*(L_0^*\subseteq L_0)$, $\beta$ the unknown vector of parameters of interest, for example, in equation \ref{eq5}, $W(a,a^*,L_0^*)=(1,a,a^*,a.a^*)^{\mathsf{T}}$.
Instead, when one is interested in effect modification by the confounder ($L_0$), then equation \ref{eq5} can be modified to \cite{lg12},
\begin{equation}\label{eq6}
E(Y_{a,M_{a^*}}|L_0)=\beta_0+\beta_{1} a+\beta_{2} a^*+\beta_{3} a . a^*+\beta_4 L_0+\beta_5 a.L_0+\beta_6 a^* L_0.
\end{equation}
The parameters indexing equation \ref{eq6} represent aspects of conditional direct and indirect effect. Effect modification of the direct and indirect effect by the confounder $L_0$ is captured by $\beta_5$ and $\beta_6$, respectively.
\\
In applied mediation analyses, some of the mediator-outcome confounders are often affected by the exposure. This generates complex forms of confounding, which are essentially guaranteed to exist in longitudinal mediation analyses (see the next section). In particular, this induces a violation of identification assumptions 2 and 4, thereby rendering the natural direct and indirect effect non-identified without making strong untestable assumptions \cite{rv15}. In such cases, we will focus on the identification of specific path-specific effects (for more details refer to \cite{VVR}). In particular, we redefine the counterfactual outcome $Y_{a,M_{a^*}}$ as expressing the outcome that would be observed if the exposure $A$ was set to $a$ (so that in particular $L_1$ takes on the value $L_{1a}$), and if $M$ was set to the value it would take if the exposure $A$ was set to $a^*$ and $L_1$ was still kept at $L_{1a}$; this counterfactual is more carefully denoted $Y_{a,L_{1a}, M_{a^*L_{1a}}}$, a notation which we will avoid to keep it manageable. In that case, the total effect can still be decomposed, but the direct effect $\beta_1$  in equation \ref{eq5} now corresponds to the combination of the paths $A\rightarrow L_1\rightarrow Y$, $A\rightarrow Y$ and $A\rightarrow L_1\rightarrow M\rightarrow Y$ (shown as $I$ in Figure \ref{f1} right), and the indirect effect, $\beta_2+\beta_3$ corresponds to the path $A\rightarrow M\rightarrow Y$ (shown as $II$ in Figure \ref{f1} right) \cite{VVR}.

\section{Extension to point exposure and time-varying mediator}
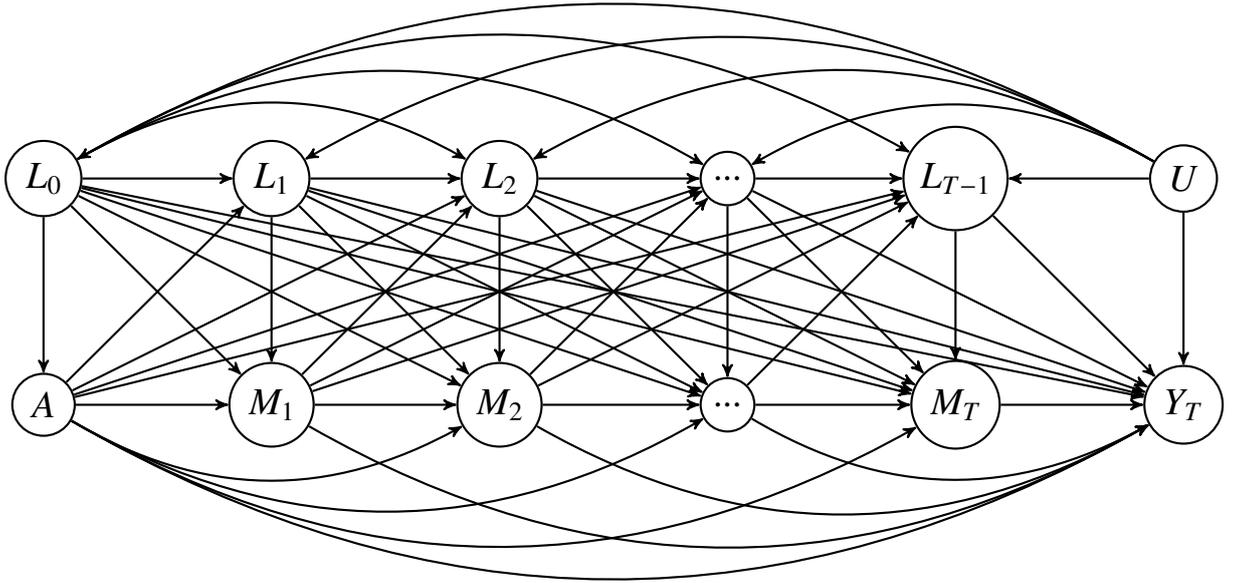
\begin{figure}
\centering
\captionsetup{justification=centering}
\begin{tikzpicture}[->,>=stealth',auto,node distance=3cm,
  thick,main node/.style={circle,draw,font=\sffamily\Large\bfseries}]

  \node[main node] (1) {$A$};
  \node[main node] (2) [above of =1]{$L_0$};
  \node[main node] (3) [right of =2]{$L_{1}$};
  \node[main node] (4) [below of=3] {$M_{1}$};
  \node[main node] (5) [right of=3] {$L_{2}$};
  \node[main node] (6) [right of=5] {$...$};
  \node[main node] (7) [right of=6] {$L_{T-1}$};
  \node[main node] (8) [below of=5] {$M_{2}$};
  \node[main node] (9) [below of=6] {$...$};
    \node[main node] (10) [below of=7] {$M_{T}$};
  \node[main node] (11) [right of=10] {$Y_{T}$};
  \node[main node] (12) [above of=11] {$U$};

  \path[every node/.style={font=\sffamily\small}]
    (2) edge node [below] {} (1)
    (2) edge node [right] {} (3)
    (2) edge node [below] {} (4)
    (2) edge node [below] {} (8)
    (2) edge node [below] {} (9)
    (2) edge node [right] {} (10)
    (2) edge node [right] {} (11)
    (1) edge node [above] {} (3)
    (1) edge node [above] {} (4)
    (1) edge node [above] {} (5)
    (1) edge node [above] {} (6)
    (1) edge node [above] {} (7)
    (3) edge node [below] {} (4)
    (3) edge node [right] {} (5)
    (3) edge node [below] {} (8)
    (3) edge node [right] {} (9)
    (3) edge node [below] {} (10)
    (3) edge node [right] {} (11)
    (4) edge node [right] {} (5)
    (4) edge node [right] {} (6)
    (4) edge node [right] {} (7)
    (4) edge node [right] {} (8)
    (5) edge node [below] {} (6)
    (5) edge node [right] {} (8)
    (5) edge node [below] {} (9)
    (5) edge node [right] {} (10)
    (5) edge node [below] {} (11)
    (6) edge node [right] {} (9)
    (6) edge node [right] {} (7)
    (6) edge node [right] {} (10)
    (6) edge node [right] {} (11)
    (7) edge node [right] {} (10)
    (7) edge node [right] {} (11)
    (8) edge node [right] {} (9)
    (8) edge node [above] {} (6)
    (8) edge node [above] {} (7)
    (9) edge node [right] {} (10)
    (9) edge node [above] {} (7)
    (12) edge node [left] {} (7)
    (12) edge node [below] {} (11)
   (10) edge node [right] {} (11)
   (2) edge[bend left] node [right] {} (6)
   (2) edge[bend left] node [right] {} (7)
   (4) edge[bend right] node [right] {} (11)
   (8) edge[bend right] node [right] {} (11)
   (9) edge[bend right] node [right] {} (11)
   (1) edge[bend right] node [right] {} (8)
   (1) edge  [bend right] node [right] {} (9)
   (1) edge[bend right] node [right] {} (10)
   (1) edge  [bend right] node [right] {} (11)
    (12) edge [bend right] node [right] {} (6)
    (12) edge [bend right] node [right] {} (2)
    (12) edge [bend right] node [right] {} (3)
    (12) edge [bend right] node [right] {} (5)
    (2) edge [bend left] node [right] {} (5);

     \end{tikzpicture}
\caption{\small{Causal Directed Acyclic Graph with exposure induced mediator-outcome confounders. There can be unmeasured confounders between any selected time periods. However, to keep the DAG simple we display the unmeasured confounder at only one time point. The vectors $L_t$ at time $t=1,\dots, (T-1)$ include $Y_t$}.}
\label{f2}
\end{figure}

In this paper, we will extend the above framework to longitudinal mediators and outcomes, while still focussing on point exposures. Let $T$ be the number of visits (excluding the baseline visit t = 0) on which individuals are measured, $t=1,2,3,...,T$. In particular, suppose that at each of these times we observe data on an outcome $Y_t$, mediators $M_t$ and a vector of covariates $L_t$. Figure \ref{f2} depicts how they are related over time. Throughout, we will denote the history of measurements up to time $t$ using a bar (e.g. $\overline{M}_t=(M_1,M_2,...,M_T)$). We will further assume that $Y_t$ may have been influenced by $\overline{M}_t$ and $\overline{L}_t$, and that $M_t$ may have been influenced by $\overline{L}_t$. We are now ready to define the nested counterfactual outcome $Y_{t,a,\overline{M}_{t,a^{*}}}$ at each time point $t$ as the outcome that would be observed at time $t$ if, for $s=1,\dots,T-1$, $L_s$ was set to the value we would have observed if the exposure was set to $a$ (and all previous instances of $L$ and $M$ to the values we have already set) and $M_s$ was set to the value we would have observed if the exposure was set to $a^*$ (and all previous instances of $L$ and $M$ to the values we have already set).
Correspondingly, the contrast $E[Y_{t,a,\overline{M}_{t,a^{*}}}-Y_{t,a^*,\overline{M}_{t,a^{*}}}]$ encodes the path-specific effects on the outcome at time $t$ of changing the exposure from $a$ to $a^*$, which captures the effect of $A$ on $Y$ along all pathways, except those where $A$ directly influences one of the mediators $M_s,s=1,...,T$. It thus captures the effect of exposure on outcome along any of the pathways in Figure \ref{f2} where $A$ directly influences either the outcome, or one of the time-varying confounders $L_t$, which may then subsequently affect outcome. The contrast $E[Y_{t,a,\overline{M}_{t,a}}-Y_{t,a,\overline{M}_{t,a^{*}}}]$ encodes the path-specific effects on the outcome at time $t$ of changing the exposure from $a$ to $a^*$, as a result of $A$ directly influencing one of the mediators $M_s$. It thus captures the effect of exposure on outcome along any of the pathways in Figure \ref{f2} where $A$ directly influences one of the time-varying mediators $M_t$, which may then subsequently affect outcome. Note that the mediated effect on which we focus, thus excludes pathways whereby treatment initially influences time-dependent patient characteristics $L$, which then in turn influence the mediator and thereby the outcome. Those pathways will be attributed to the indirect effect via those patient characteristics. This seems logical from an interpretational point of view, but is also a more fundamental requirement: the effect of treatment transmitted along the combination of all pathways that intercept one or multiple mediators (regardless of where in the causal chain it intercepts these variables) cannot be identified without making overly stringent assumptions \cite{asp05,sh13,SV17}.

Below we will discuss identification and estimation of the above path-specific effects under the no-unmeasured confounding assumptions implicit in the causal diagram of Figure \ref{f2}, which is assumed to represent a nonparametric structural equation model with independent errors. In particular as in \cite{SV17}, we will assume that the same set of baseline covariates $L_0$ is sufficient to control for confounding of associations between exposure and outcome, and of exposure and mediator. We will moreover assume that at each time $t$ and each time $s\leq t$, adjustment for the history $A,\overline{L}_s,\overline{M}_{s-1}$ suffices to control for confounding of the effect of $M_s$ on $Y_t$. Throughout, we will assume that there is no loss to follow-up and measurement error in the mediators.

\subsection{Natural Effect Models for Longitudinal Mediators}
Previously reviewed natural effect models can be extended to longitudinal mediators and outcomes. For instance, equation \ref{eq5} can be extended to accommodate time dependent mediators as follows:
\begin{equation}\label{eq7}
g\{E[Y_{t,a,\overline{M}_{t,a^*}}]\}=\alpha_0+\alpha_{1}a+\alpha_{2}a^*+\alpha_{3}t+\alpha_{4}t.a+\alpha_{5}t.a^*
\end{equation}
for a user-specified link function $g(.)$.
This generalises the marginal structural model
\[
g\{E[Y_{t,a}]\}=\alpha_0+(\alpha_{1}+\alpha_{2})a+\alpha_{3}t+(\alpha_{4}+\alpha_{5})t.a
\]
When $g(.)$ is the identity link, then according to equation \ref{eq7} the total effect of a unit increase in the exposure, $\alpha_1+\alpha_2+(\alpha_4+\alpha_5)t$, can be decomposed into direct effect
\[
E\{Y_{t,1,\overline{M}_{t,0}}-Y_{t,0,\overline{M}_{t,0}}\}=\alpha_1+\alpha_{4}t
\]
and the indirect effect
\[
E\{Y_{t,1,\overline{M}_{t,1}}-Y_{t,1,\overline{M}_{t,0}}\}=\alpha_2+\alpha_{5}t.
\]
Equation \ref{eq7} excludes the possibility of mediator-exposure interaction (on the scale of the link function $g(.)$). However, such interactions can be allowed in the model by including $a.a^*$ term in the model.

\subsection{Estimation}
Here we describe how the coefficients in natural effect models of the form (equation \ref{eq7}) can be estimated for a binary point exposure, a longitudinal outcome, and a longitudinal categorical mediator, using inverse probability weighting. Estimation using natural effect models is not limited to these types of variables but can in principle accommodate arbitrary exposures, mediator and outcomes  (e.g. binary, categorical, continuous).
\subsubsection{Regression for Exposure}
To adjust for confounding of the exposure-mediator and exposure-outcome associations, we will first calculate inverse probability of exposure weights. Assuming
\begin{eqnarray}\label{eq8}
\textit{logit}\left( P\left[A=1|L_0=l_0\right]\right)&=& \gamma_0+\gamma_1 l_0,
\end{eqnarray}
estimated probabilities ($\hat{p}_i$) for each individual $i$ can be obtained as $\hat{p}_i=\text{expit}(\hat{\gamma}_0 +\hat{\gamma}_1 l_{0i})$. The weight for the $i^{th}$ individual is then $w_{i}^{a}=\frac{1}{\hat{p}_i}$ if $A_i=1$ and $w_{i}^{a}=\frac{1}{(1-\hat{p}_i)}$ if $A_i=0$.
\subsubsection{Regression for the mediators}
To decompose the total effect into path-specific effects, while adjusting for confounding of the mediator-outcome associations, we will next calculate mediator weights. For a categorical mediator $M_t$ with possible values $k= 0,...,K$, we will fit multinomial models at each time point. For the mediator at the first time point, we may for instance fit the multinomial logistic model:
\begin{eqnarray}\label{eq9}
pr[M_{1}=k|A=a,\overline{L}_1=\overline{l}_1]&=& \frac{\exp(\delta_{0k}+\delta_{1k} a+\delta_{2k} \overline{l}_1)}{1+\sum_{l=1}^{K} \exp(\delta_{0l}+\delta_{1l} a+\delta_{2l} \overline{l}_1)}, k=1,2,...,K \nonumber \\
pr[M_1=0|A=a,\overline{L}_1=\overline{l}_1]&=& \frac{1}{1+\sum_{l=1}^{K} \exp(\delta_{0l}+\delta_{1l} a+\delta_{2l} \overline{l}_1)}
\end{eqnarray}
and then subsequently for all later time points:
\begin{eqnarray}\label{eq10}
pr[M_{t}=k|A=a,\overline{M}_{t-1}=\overline{m}_{t-1},\overline{L}_{t}=\overline{l}_t]&=& \frac{\exp(\delta_{0k}+\delta_{1k} a+\delta_{2k} \overline{m}_{t-1}+\delta_{3k} \overline{l}_t)}{1+\sum_{l=1}^{K} \exp(\delta_{0l}+\delta_{1l} a+\delta_{2l} \overline{m}_{t-1}+\delta_{3l} \overline{l}_t)}, k = 1,2,...,K \nonumber \\
pr[M_t=0|A=a,\overline{M}_{t-1}=\overline{m}_{t-1},\overline{L}_t=\overline{l}_t]&=& \frac{1}{1+\sum_{l=1}^{K} \exp(\delta_{0l}+\delta_{1l} a+\delta_{2l} \overline{m}_{t-1}+\delta_{3l} \overline{l}_t)}
\end{eqnarray}
The mediator weight (Appendix 1) for the $i^{th}$ individual at time $t$ corresponding to counterfactual $Y_{ta\overline{M}_{ta^*}}$ is then
\begin{equation}\label{eq11}
w_{t,i}^m (a^*)=\prod_{s=1}^{t}\frac{P(M_{t}=m_{s,i}|A=a^{*},\overline{L}_{s}=\overline{l}_{s,i},\overline{M}_{s-1}=\overline{m}_{s-1,i})}{P(M_{t}=m_{s,i}|A=a,\overline{L}_{s}=\overline{l}_{s,i},\overline{M}_{s-1}=\overline{m}_{s-1,i})}.
\end{equation}
The validity of this approach follows from references \cite{zv18} and \cite{SV17}. While reference \cite{zv18} focuses on interventional analogues to direct and indirect effect, it provides identical identification results as in \cite{SV17}, rendering the inverse weighting strategy in \cite{zv18} applicable. Note that \cite{zv18} does not consider estimation of natural effect models, however.
\subsubsection{Fitting the natural effect models}
For binary exposure (coded 0 or 1), the general process for estimating the natural effect models can be described as follows:
\begin{enumerate}
\item Fit a suitable model for the exposure (e.g. equation \ref{eq8}) on the original data.
\item For each time $t=1,...,T$, fit a suitable model for the mediator (e.g. if $M_t$ is categorical, one may use equation \ref{eq10}) conditional on the history of previous mediators, exposures and time-varying confounders, where the latter may include the outcome history. This can also be based on one model across all times fitted using pooled regression.
\item Construct a new data set by replicating each observation in the original data set twice for each time point and include an additional variable $a^*$, where $a^*$ is equal to the original exposure ($A$) for the first replication and equal to $1-A$ for the second replication. In addition add an identification variable to indicate which data rows originate from the same subject, along with a visit time indicator, $t=1,2,...,T$ (see Table 1).
 \begin{table}
\caption{Schematic display of the weighting in longitudinal mediation analysis}\label{t-1}
\begin{center}
\begin{tabular}{ccccccc}
%p{1cm}|p{1cm}|p{1cm}|p{1cm}|p{1cm}|p{1cm}|p{1cm}
 \hline
 Id & $A_i$ & \boldsymbol{${a}$} & $a^*$ &T& $Y_i$& $w_i$ \\
 \hline
1 & 1 & 1 & 1&$t_1$&$Y_{11}$& $\frac{1}{\hat{p}_{1}}$ \\
1 & 1 & 1 & 0& $t_1$&$Y_{11}$& $\frac{w_{1,i}^m(a^*)}{\hat{p}_{1}}$\\
1 & 1 & 1 & 1&$t_2$&$Y_{12}$& $\frac{1}{\hat{p}_{1}}$\\
1 & 1 & 1 & 0& $t_2$&$Y_{12}$& $\frac{w_{2,i}^m (a^*) w_{1,i}^m (a^*)}{\hat{p}_1}$\\
\vdots&\vdots&\vdots&\vdots&\vdots&\vdots&\vdots\\
2 & 0 & 0 & 0&$t_1$&$Y_{21}$& $\frac{1}{1-\hat{p}_2}$\\
2 & 0 & 0 & 1& $t_1$&$Y_{21}$& $\frac{w_{1,i}^m (a^*)}{1-\hat{p}_2}$\\
2 & 0 & 0 & 0&$t_2$&$Y_{22}$& $\frac{1}{1-\hat{p}_2}$\\
2 & 0 & 0 & 1& $t_2$&$Y_{22}$& $\frac{w_{2,i}^m(a^*) w_{1,i}^m(a^*)}{1-\hat{p}_2}$\\
 \hline
\end{tabular}
\end{center}
\end{table}
\item At time $t$ the weight for the $i^{th}$ individual corresponding to entry $a^*=0, 1$ in Table \ref{t-1} is computed as
\[w_{t,i} (a^*)=w_i^a w_{t,i}^m (a^*).\]
\item Fit the natural effect models (equation \ref{eq7}) by regressing the outcome on time $t$, $a$ and $a^*$ and the time interactions ($ta, ta^*$) on the basis of expanded data set, using weighted Generalised Estimating Equations (GEE) with independence working correlation. The weights in this regression are the weights ($w_{t,i}(a^*)$) computed in the previous step. Note that the use of an independence working correlation is critical to ensure that the right weights are assigned to the right records \cite{sv07}.
\end{enumerate}
The above procedure can easily be implemented in standard software (Appendix 2). When the natural effect models include baseline covariates, e.g. $g\{E[Y_{t,a,\overline{M}_{t,a^*}}|L_0]\}=\alpha_0+\alpha_1 a+\alpha_2a^*+\alpha_3t+\alpha_4 L_0+\alpha_5t.a+\alpha_6 t.a^*+\alpha_7 a.L_0 +\alpha_8 a^{*}.L_0$, then $w_i^a$ can be set to 1 for all individuals because the adjustment for confounding by $L_0$ now happens via a standard regression adjustment \cite{lg12,sv12}. The resulting weights, $w_{t,i}(a^*)$ will typically be more stable \cite{lg12,sv12}.
\section{Estimation of standard errors and confidence intervals}
To compute confidence intervals for the parameters indexing a natural effect models, it is tempting to rely on the standard GEE output. However, unlike with weighting procedures for marginal structural models, this does not guarantee conservative intervals. The bootstrap forms an attractive alternative, but may have the drawback that when the mediators are categorical, some of the categories might not be selected in every bootstrap replication. It is for this reason we give an alternative method for computing the variance \cite{DHY}, which forms a hybrid between the use of robust standard error estimators from the GEE output, and the parametric bootstrap. The proposed procedure for computing the correct confidence intervals and standard error is described as follows:
\begin{itemize}
\item[Step-1] Fit a model for the exposure (e.g. logistic if $A$ is binary) and extract the regression coefficients and their covariance matrix.
\item[Step-2] Fit a model for the mediator (e.g multinomial if $M$ is categorical), for each time point as described in the above section, and extract the regression coefficients and their covariance matrix.\\
\end{itemize}
Repeat the following procedure $B$ (e.g. $B=1000$) times:
\begin{itemize}
\item[Step-3] Randomly perturb the estimated coefficients from steps 1 and 2 by adding mean zero normal noise with covariance matrix as estimated in steps 1 and 2.
\item[Step-4] Using the sampled coefficients from step 3, recompute the weights as described in equation \ref{eq11}.
\item[Step-5] Using these new set of weights, estimate the coefficients corresponding to $a$, $a^*$, $ta$, and $ta^*$ in the natural effect models.
\item[Step-6] In the $j^{th}$ iteration, for $j=1,...,B$, store the regression coefficients, $\hat{\alpha}^{(j)}$, and their variance-covariance, $Var(\hat{\alpha}|\hat{w}_j)$, from step 5. Here, the notation $Var(\hat{\alpha}|\hat{w}_j)$ expresses the uncertainty at fixed weights.
\item[Step-7] Using the results from these $B$ iterations estimate the variance-covariance of the regression coefficients corresponding to $a$, $a^*$, $ta$, and $ta^*$ using the formula
\[\frac{1}{B}\sum_{j=1}^B Var(\hat{\alpha}|\hat{w}_j)+\frac{1}{B-1}\sum_{j=1}^B \left(\hat{\alpha}^{(j)} - B^{-1}\sum_{k=1}^B \hat{\alpha}^{(k)}\right)\left(\hat{\alpha}^{(j)} - B^{-1}\sum_{k=1}^B \hat{\alpha}^{(k)}\right)'.\]
\item[Step-8] A confidence interval for the $k^{th}$ component of $\alpha$ can be obtained as $\hat{\alpha}_k\pm1.96\times se_k$, where the standard error $se_{k}$ is the square root of the element in the $k^{th}$ row and $k^{th}$ column of the matrix computed in step 7.
\end{itemize}
\section{Illustrating Example}
As an illustrating example we use the data from the Millennium Cohort Study (MCS), a longitudinal study of children born in the UK between September 2000 and January 2002, which has been described elsewhere \cite{cp14}. Ethical approval for the MCS was received from a Research Ethics Committee at each time point of the survey round. Data were obtained from the UK Data Archive, University of Essex in March 2014. The first study contact with the cohort child was at 9 months, with survey interviews carried out by trained interviewers in the home with the main respondent (usually the mother) and their partner. Information was collected on 72$\%$ of those approached, providing information on 18,818 infants. For the analysis, data was used from four waves, when the children were aged 3, 5, 7 and 11 years. The total number of respondents who had participated in all four waves (by the age of 11 years of child) had declined to 10,313 (56$\%$ of the respondents who had taken part at the 9 months time period). To illustrate the techniques developed in this paper we use complete cases, $N= 5188$. Basic descriptives for some of the variables considered (e.g. Strengths and Difficulty Questionnaire at 3, 5, 7 and 11) show that the distribution of ``Yes" categories in complete cases are slightly under-represented compared to the distributions of ``Yes" in observed cases \cite{hm18}; similar is the observation for the maternal psychological distress variable. Further, for baseline confounders such as maternal education, missing cases were primarily from lower education categories such as ``O-Levels/GCSE" and ``Others".

In our illustration, we are interested in examining the extent to which socio-economic position (SEP) in infancy and maternal psychological distress measured over the period of child development (ages 3 5 7 and 11 years) affect child mental health measured at the age of 11 years. The baseline confounders considered for SEP-child mental health relation were maternal education ($L01$), ethnicity ($L02$), gender ($L03$), marital status ($L04$) and age of mother ($L05$). Even though this is not a very comprehensive list of baseline confounders, we believe it to be a predominant subset of confounders of the relation between SEP and child mental health. The exposure (SEP) was categorised as ``Above 60$\%$ median'' and ``below 60$\%$ median''. Maternal psychological distress was measured using the Kessler-6 scale, a self-reported measure collected from the mother of the child by asking how often in the past 30 days she had felt:  ``So depressed that nothing could cheer you up", ``Hopeless", ``Restless or Fidgety", ``That everything was an effort", ``Worthless" and ``Nervous". Each of these items had a five-point response from  ``None of the time" (0) to ``All of the time" (4). Responses to each item were combined to produce a single score ranging from 0 to 24. Scores exceeding 8 were truncated at 8. Child mental health problems were assessed using the Strengths and Difficulties Questionnaire, a 25-item measure completed by the mother. The overall score was classified into two groups  ``Normal", 0, if the child had a score in between 0 and 13 and ``Border-line abnormal", 1, if they had score in between 14 and 40. These cut-offs were chosen based on recommended guidelines \cite{gr97}.

In this illustration we present two scenarios where we take the outcome measured at the end of study, and the outcome measured from each time period. First, we estimate the effect of the parental SEP when the child's age was 9 months, on the child's mental health at age 11, with the mediators being  maternal psychological distress measured at birth.  The time-varying confounders include the time-varying outcomes and mediators measured at ages 3, 5 and 7, which we believe to be the main confounders, although we acknowledge that there may be other confounding factors such as family environment, medical treatment and others, which were unavailable. We then fitted a logistic natural effect models without the interaction term between exposure and mediator (P = 0.968).  The indirect effect $\beta_2$ in this model captures the combined effect along all pathways whereby $A$ directly affects one of the mediators, which may then subsequently affect outcome at age 11. The direct effect $\beta_1$ captures the combined effect along all remaining pathways.

To estimate the direct and indirect effect on the end of study outcome, the data needs to be arranged as shown in Table \ref{t-2}. Results from this analysis show that SEP in childhood has a bigger effect compared to the effect via maternal psychological distress. The effect through SEP is 0.544 (corresponding to an odds ratio of exp(0.544)=1.72, 95$\%$ confidence interval (CI) 1.27 to 2.33). The effect via maternal distress is 0.122 (corresponding to an odds ratio of exp(0.122)=1.13, 95$\%$ CI 1.07 to 1.20), thus indicating that most of the effect is not via maternal distress.
 \begin{table}
\caption{Schematic display of the weighting for the end of study outcome}\label{t-2}
\begin{center}
\begin{tabular}{cccccc}
%p{1cm}|p{1cm}|p{1cm}|p{1cm}|p{1cm}|p{1cm}|p{1cm}
 \hline
 Id & $A_i$ & $a$ & $a^*$ & $Y_i$& $w_i$ \\
 \hline
1 & 1 & 1 & 1&$Y_{14}$& $\frac{1}{\hat{p_{1}}}$ \\
1 & 1 & 1 & 0&$Y_{14}$& $w_{i}(a^*)$\\
2 & 0 & 0 & 0&$Y_{24}$& $\frac{1}{1-\hat{p_{2}}}$ \\
2 & 0 & 0 & 1&$Y_{24}$& $w_{i}(a^*)$\\
\vdots&\vdots&\vdots&\vdots&\vdots&\vdots\\
 \hline
\end{tabular}
\end{center}
\end{table}

\begin{figure}
\includegraphics[height=0.5\textheight,width=1\linewidth]{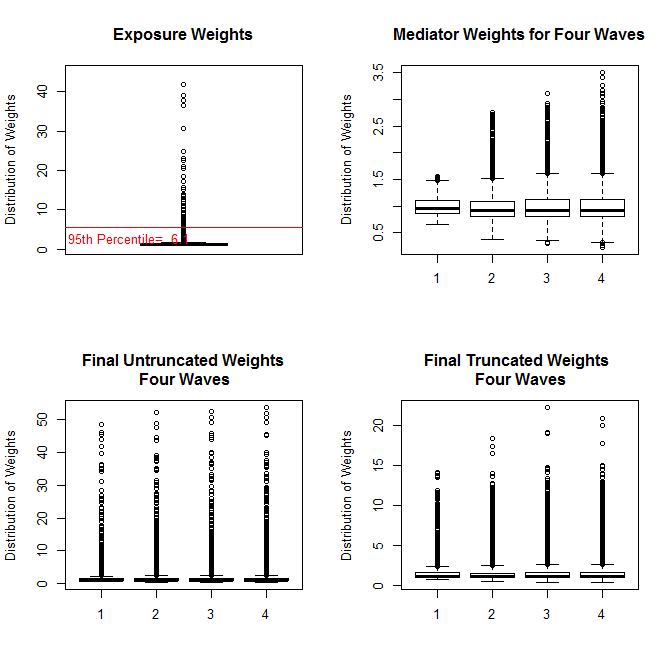}
\caption{Distribution of exposure weights (upper left panel),  mediator weights (upper right panel), untruncated weights (lower left panel) and the final 95\% truncated weights (lower right panel).}
\label{f3}
\end{figure}

We next estimated the direct and indirect effect on all outcomes. For this purpose, the data needs to be arranged as in Table \ref{t-1}.
Using logistic regression for the binary outcome, we estimated the time-specific effects.  The distribution of weights for each time point is given in Figure \ref{f3}.

%\begin{figure}
%\includegraphics[height=0.5\textheight,width=1\linewidth]{mweights.png}
%\caption{Distribution of mediator weights for four waves}
%\label{f3}
%\end{figure}

%\begin{figure}
%\includegraphics[height=0.5\textheight,width=1\linewidth]{mfw.png}
%\caption{Distribution of final untruncated and truncated 95 weights for four waves}
%\label{f4}
%\end{figure}

\begin{figure}
\centering
\includegraphics[height=0.5\textheight,width=.8\linewidth]{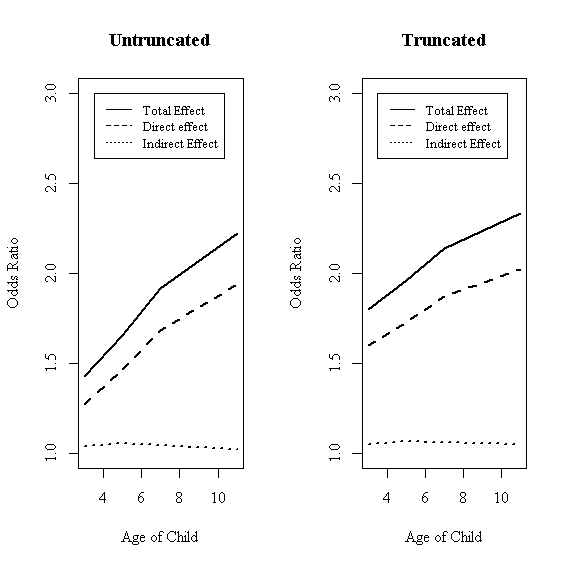}
\caption{Trend in direct, indirect and total effect of socio-economic position and maternal psychological distress on child mental health status, using untruncated (left) and truncated weights (right).}
\label{f4}
\end{figure}

\begin{figure}
\begin{subfigure}[a]{0.8\textwidth}
\centering
\includegraphics[width=.8\linewidth]{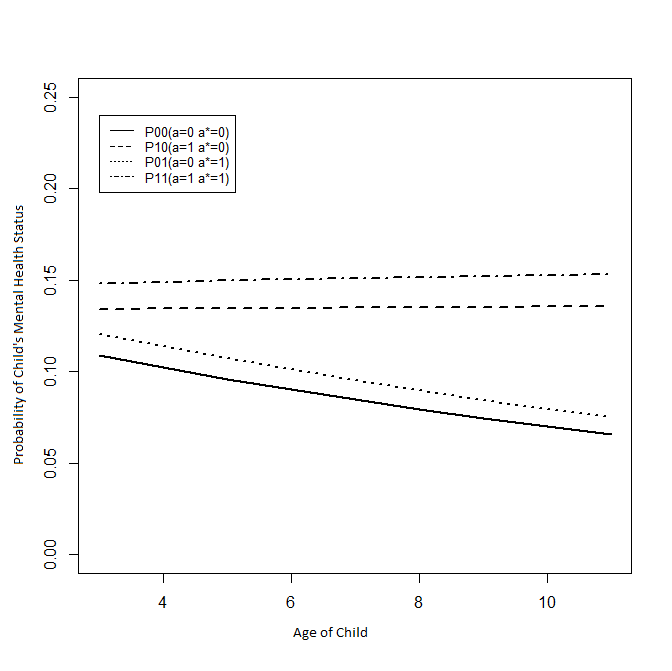}
%\caption{Graphical depiction of Probability of child's mental health status, using untruncated weights}
%\label{f5a}
\end{subfigure}
\begin{subfigure}[b]{0.8\textwidth}
\centering
\includegraphics[width=.8\linewidth]{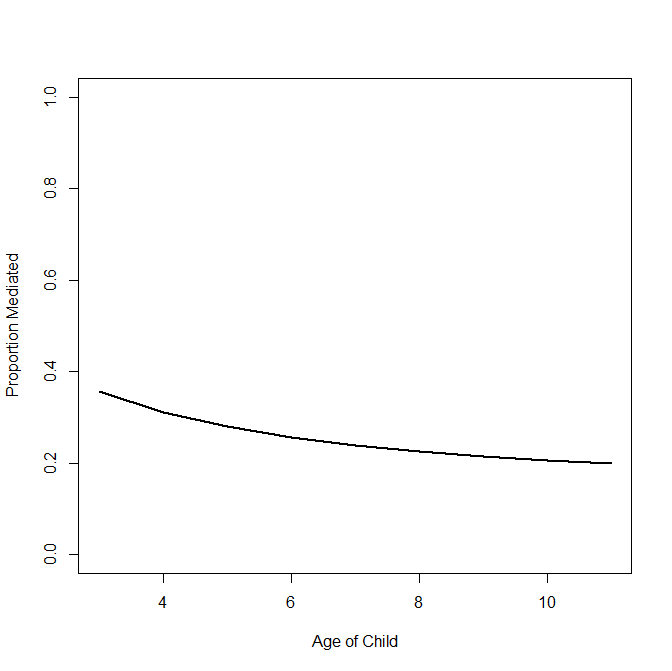}
%\caption{Graphical depiction of proportion mediated using untruncated weights}
%\label{f5b}
\end{subfigure}
\caption{Graphical presentation of the probability of child's mental health status (top) and proportion mediated using untruncated weights (bottom).}
\label{f5}
\end{figure}

The estimates of the time/path specific effects are presented in Table \ref{t4}.
\begin{table}
\centering
\caption{Odds ratio estimates from longitudinal mediation analysis using natural effect models}\label{t4}
\begin{center}
\begin{tabular}{ |p{2cm}|p{1cm}|p{1.75cm}|p{1cm}|p{1.75cm}|p{1cm}|p{1.75cm}|p{1cm}|p{1.75cm}|}
 \hline
 \multicolumn{9}{|c|}{Effect (95\% confidence interval (CI))} \\
 \hline
 Estimate& \multicolumn{2}{c}{Age 3}&\multicolumn{2}{|c}{Age 5} & \multicolumn{2}{|c}{Age 7} &\multicolumn{2}{|c|}{Age 11}\\
 \hline
\multicolumn{9}{|c|}{Untruncated weights} \\
 \hline
 &Odds ratio& 95$\%$ CI&Odds ratio&95$\%$ CI&Odds ratio&95$\%$ CI&Odds ratio&95$\%$ CI\\
 \hline
 Direct &1.27&0.90, 1.80&1.47&1.10, 1.95&1.68&1.28, 2.21&1.94&1.41, 2.67\\
 ($\alpha_1+\alpha_4\times t$)&&&&&&&&\\
 Indirect &1.12&1.04, 1.21&1.13&1.06, 1.21&1.14&1.06, 1.23&1.15&1.02, 1.28\\
 ($\alpha_2+\alpha_5\times t$)&&&&&&&&\\
 Total&1.43&1.01, 1.79& 1.66 &1.24, 1.95&1.92&1.46, 2.21&2.22&1.63, 2.67\\
 \hline
 \multicolumn{9}{|c|}{Truncated 95$\%$ weights}\\
 \hline
 Direct&1.59&1.28, 1.99&1.73&1.44, 2.08&1.87&1.54, 2.27&2.02&1.58, 2.60\\
 ($\alpha_1+\alpha_4\times t$)&&&&&&&&\\
 Indirect&1.13 &1.05, 1.21&1.14&1.06, 1.21&1.15&1.05, 1.23&1.15&1.05, 1.27\\
 ($\alpha_2+\alpha_5\times t$)&&&&&&&&\\
 Total&1.80&1.45, 2.23&1.96 &1.65, 2.33&2.14&1.79, 2.55&2.33&1.86, 2.92\\
 \hline
\end{tabular}
\end{center}
\end{table}
From Table \ref{t4} we note that both the maternal mental health and SEP contributed to the well-being of a child. Table \ref{t4} provides the odds ratios from the fitted natural effect model. For instance, for an individual child, altering the level of socio-economic position from high (=0) to low (=1), while controlling the maternal psychological distress at the levels as naturally observed at any given level of SEP, $a$, almost doubles (exp$(\alpha_1+\alpha_4*3=0.662)=1.94$; 95 $\%$ CI  1.41, 2.67) the odds of displaying mental health issues at age 11. This is much more pronounced than at age 3 (exp$(\alpha_1+\alpha_4*0=0.241)=1.27$; 95 $\%$ CI 0.90, 1.80). This may be partly explained by the fact that the mental health at age 3 was reported by parents and not the children themselves. Similarly, altering levels of maternal psychological distress as observed at low SEP to levels that would have been observed at high maternal distress while controlling their SEP at any given level, $a$, increases the odds of a child suffering from mental health issues at age 11 (exp$(\alpha_2+\alpha_5*3=0.135)=1.14$; 95$\%$ CI: 1.02, 1.28). Similar results were obtained using truncated weights. The truncated weights were computed by resetting the value of weights greater than $95^{th}$ percentile (weight = 6.1) to the value of the $95^{th}$ percentile (6.1). Figure \ref{f4} shows the trend in the direct and indirect effect over the period. It can be seen that the total effect\ (SEP and history of maternal psychological distress) increases over the period. In addition, the direct effect of SEP on child mental health was observed to have long term effects. On the other hand, the effect via maternal psychological distress on child mental health was almost constant. To estimate the proportion mediated we first computed the probabilities, $P_{00}, P_{01}, P_{10}$, and $P_{11}$, where $P_{kl}$ is the probability of boarder line abnormal child mental health when $a$ is set to $k = 0, 1$ and $a^*$ is set to $l = 0, 1$  calculated using equation \ref{eq7} with logit link. For example $P_{a.a^*}$ was computed as $=\frac{exp(\alpha_0+\alpha_1 a+\alpha_2 a^*+\alpha_3 t+\alpha_4 ta+\alpha_5ta^*)}{1+exp(\alpha_0+\alpha_1 a+\alpha_2 a^*+\alpha_3 t+\alpha_4 ta+\alpha_5ta^*)}$. Figure \ref{f5} (top panel)  shows the plot of these probabilities and the proportion mediated $\frac{(P_{11}-P_{10})}{(P_{11}-P_{00})}$, using untruncated weights. From Figure \ref{f5} (bottom panel) we note that the proportion mediated declines over time. Inference does not seem to change when using truncated weights (see Figure \ref{f4}, right panel).

\section{Discussion}
We have generalised the popular class of marginal structural models, which parameterise the effect of a point exposure on a (possibly time-varying) outcome, to enable effect decomposition w.r.t. a time-varying mediator in the presence of time-varying confounding. We have focussed on the situation where the mediator was multinomial, although the proposed inverse weighting strategy for fitting these models is readily adapted to mediators of arbitrary type (e.g. binary, continuous, count). Caution is warranted, however, when working with continuous mediators, because the proposed approach then requires inverse weighting by the mediator density, which can easily lead to instability. While we have focussed on dichotomous exposures, the proposal is also readily extended to arbitrary exposure distributions. As in Vansteelandt et al.\cite{sv12}, this requires augmenting the observed data set by picking for example, five random values $a^*$ from the exposure distribution and next following steps 2-6 described in the above section. To avoid inverse weighting by the exposure density, one may then consider standard regression adjustment for confounding of the exposure-outcome association.

The validity of the proposed approach critically relies on the assumptions expressed by Figure 2. In particular, we assume that adjustment for the measured confounders $L_0$ suffices to identify the effect of exposure on all mediators, time-varying confounders and outcome. We moreover assume that no unmeasured variables are simultaneously associated with the mediator on the one hand, and time-varying confounders or outcome on the other hand. We thus in particular assume that the effect of time-varying confounders at time $t$ on future mediators is unconfounded after adjusting for the history of all available measurements prior to time $t$. More importantly, we assume that the effect of mediators at time $t$ on future time-varying confounders and outcome is unconfounded after adjusting for the history of all available measurements prior to time $t$. These assumptions are strong and may be especially implausible when the available number of time-varying confounders is large, as it may then be difficult to believe that all of them are associated with future mediators only by means of a causal effect. Similar assumptions can for instance be avoided when the aim is to assess the overall effect of the time-varying mediators on outcome. Traditional approaches for longitudinal mediation analysis (e.g. MacKinnon, 2008 \cite{mk08}) also necessitate such assumptions, although they may not explicate them; in fact,  many approaches assume in addition that also the effect of time-varying confounders at time $t$ on future outcomes is unconfounded after adjusting for the history of all available measurements prior to time $t$, which is especially difficult to believe when there is confounding by past outcomes.

Like other approaches for path-specific effects \cite{asp05}, the proposed approach additionally relies on so-called cross-world counterfactual independencies \cite{Robins,N15}. These can be viewed as a strengthening of the assumptions listed in the previous paragraph, in the following sense. For instance, with a single mediator and a dichotomous exposure, assumption (2) (together with assumption (1) and substituting $a$ by $a^*$) implies that $Y_{a^*m}\ci M_{a^*}|L_0$. This shows that the cross-world independence assumption (4) is closely related to the ignorability assumption (2). In particular, assumption (2) expresses that subjects who would have low versus high levels of the mediator if given exposure $a^*$ are exchangeable (within strata of $L_0$) in terms of what their outcome would be if given exposure $a^*$ and if their mediator were set to $m$. Assumption (4) expresses that such individuals should additionally be exchangeable in terms of what their outcome would be if given exposure $a$ and if their mediator were set to $m$. In our opinion, it is usually hard to imagine that investigators could have such detailed knowledge to understand that assumption (2) holds, but not (4). In particular, if the investigators judge subjects with low versus high levels of the mediator at a given exposure $a^*$ to be exchangeable (within strata of $L_0$), then we would generally believe them to be exchangeable in terms of both $Y_{a^*m}$ and $Y_{am}$. In that sense, we personally view the additional cross-world counterfactual independencies as relatively more innocent. Even so, the need for cross-world counterfactual independence assumptions is frustrating, as it implies that the aforementioned exchangeability cannot be avoided, even if data were available from multiple randomised experiments which randomised either exposure or mediators (or confounders), as it is needed to tie together the effects over the different paths.

Despite this, by interpreting the obtained (in)direct effect as interventional (in)direct effect \cite{VD17,zv18}, one can avoid the need for cross-world independence assumptions, along with the assumption that the effect of time-varying confounders at time $t$ on future mediators is unconfounded after adjusting for the history of all available measurements prior to time $t$. Interventional direct effect express the effect of changing the exposure while fixing the repeated mediators at random draws from their conditional distributions at each time, given the history of the measurements observed until that time (while fixing the exposure at a given level, the previous mediators at the previously drawn values, and the previous time-dependent confounders at the levels that would naturally arise under such intervention). These effects differ from the interventional effects in \cite{vtlm17,lin17} who consider random draws from the distribution of the mediator at a certain exposure level conditional only on baseline covariates. By conditioning only on baseline covariates, these draws are less representative for what that an individual might have ``naturally" experienced, making these estimands arguably less suitable to develop insight into mechanism \cite{SV17}.

In summary, we have described a simple, generic procedure for estimating direct and indirect effect of an exposure on an outcome in a longitudinal setting. The procedure can be applied in standard software (see Appendix 2 for an R implementation). As with all inverse weighting methods, monitoring of the distribution of the inverse probability weights is recommended to detect possible instabilities. In future work, we will extend the proposal to time-varying exposures, as well as examine the possibility of using fitting strategies for the inverse probability weights aimed at preventing instabilities.

\section*{Acknowledgement}
MNM and SV equally contributed to both development and writing of this manuscript. MNM was funded by the Endeavour Executive scholarship from the Department of Education, Australia, to carry out the work presented. We would like to thank Dr. Steven Hope from the Institute of Child Health, University College London, London, UK for sharing the UK Millennium Cohort data. The authors would like to declare no conflict of interest.

\bibliographystyle{vancouver}
\bibliography{LManalysis}

\begin{appendices}
There are three parts to this appendix 1) Detailed derivations of the weight procedure described in main text of the paper 2) R-code to conduct both the analysis using the proposed method and code for computing the variance using the proposed method in section 4 of the main paper and 3) Derivation of the variance method, code to generate data and conduct simulations to compare the performance of the perturbed bootstrap method.

\section*{Appendix 1}
This appendix presents detailed derivations of the procedure described in the main text. The derivation is closely related to the theorems in Lange et.al (10), but here we extended to longitudinal mediators for estimating natural direct and indirect effects.

We assume that the causal diagram of Figure 2 holds and represents a non-parametric structural equation model with independent errors (13).
It then follows from Vansteelandt (9)  that we can link $E[Y_{t(a,\overline{M}_{t,a^*})}]$ to the observed data as:
\noindent
\small{
\begin{eqnarray*}
E[Y_{t(a,\overline{M}_{t,a^*})}]
&=&
\sum_{\overline{m}_T}\sum_{\overline{l}_{T}}E[Y|A=a,\overline{M}_{T}= \overline{m}_T,\overline{L}_{T}=\overline{l}_T]\left\{\prod_{t=1}^{T}P[M_{t}=m_t|A=a^*,\overline{M}_{t-1}=\overline{m}_{t-1},\overline{L}_t=\overline{l}_t]\right.\\
&&\left. \times P[L_{t}=l_{t}|A=a,\overline{M}_{t-1}=\overline{m}_{t-1},\overline{L}_{t-1}=\overline{l}_{t-1}]\right\}P(L_0=l_0)\\
&&=\sum_a \sum_{\overline{m}_T}\sum_{\overline{l}_{T}}E[Y|A=a,\overline{M}_{T}= \overline{m}_T,\overline{L}_{T}=\overline{l}_T]
\left\{I(A=a)/P(A=a|L_0=l_0)\right\}P(A=a|L_0=l_0)
\\
&& \times \left\{\prod_{t=1}^{T}P[M_{t}=m_t|A=a^*,\overline{M}_{t-1}=\overline{m}_{t-1},\overline{L}_t=\overline{l}_t]P[L_{t}=l_{t}|A=a,\overline{M}_{t-1}=\overline{m}_{t-1},\overline{L}_{t-1}=\overline{l}_{t-1}]\right\}P(L_0=l_0)\\
&&=\sum_a \sum_{\overline{m}_T}\sum_{\overline{l}_{T}}E[Y|A=a,\overline{M}_{T}= \overline{m}_T,\overline{L}_{T}=\overline{l}_T]
\left\{I(A=a)/P(A=a|L_0=l_0)\right\}P(A=a|L_0=l_0)
\\
&& \times \left\{\prod_{t=1}^{T}P[M_{t}=m_t|A=a,\overline{M}_{t-1}=\overline{m}_{t-1},\overline{L}_t=\overline{l}_t]P[L_{t}=l_{t}|A=a,\overline{M}_{t-1}=\overline{m}_{t-1},\overline{L}_{t-1}=\overline{l}_{t-1}]\right\}P(L_0=l_0)\\
&& \times \left\{\prod_{t=1}^{T}\frac{P[M_{t}=m_t|A=a^*,\overline{M}_{t-1}=\overline{m}_{t-1},\overline{L}_t=\overline{l}_t]}{P[M_{t}=m_t|A=a,\overline{M}_{t-1}=\overline{m}_{t-1},\overline{L}_t=\overline{l}_t]}\right\}P(L_0=l_0)\\
&&=E[YI(A=a) W]
\end{eqnarray*}
}
where we define $M_0$ to be empty and where $W$ refers to the product of the weights $W_i^a W_{t,i}^{m}(a^*)$. This motivates the proposed weighting procedure.
For the estimators to have good performance infinite weights must be avoided. Our proposal therefore relies on the positivity assumption that $P(m_t|A=a,\overline{M}_{t-1}, \overline{L}_{t})>\sigma>0$ with probability 1 for all $t$, $m_{t}, a$ and $P(A=a|L_0)>\epsilon>0$ with probability 1 for all $a$.

\section*{Appendix 2}
\subsection*{Data structure}
Initially we load the original data set into variable named `mydata', with observed exposure (A), outcome ($Y_1,Y_2,Y_3,Y_4$), mediator measured at four time points ($M_1, M_2, M_3   \text{and} \ M_4$) and the confounders ($L_0$). Next we construct new variables id, and `Astar' where the latter corresponds to the value of the exposure relative to the indirect path. We also replicate a new data set by replicating the original two times such that `Astar' takes the two different possible values through the replications. We repeat this process for each mediator.
\begin{verbatim}
First read the data into a matrix/dataframe/vector named mydata
mydata<-read.csv(file="mentalhealth.csv",header=TRUE,sep=",")
N<-nrow(mydata)
#creating an Id variable
mydata$id<-1:N
#L0 has five variables hence L01 L02 L03 L04 L05.
# M1 M2 M3 M4 mediator observed at times 0 1 2 and 3.
head(mydata, n=10L)

Data could not be provided due to confidentiality issues.

#We have mediator measured at four time points.
#To implement the method described in the paper
#we duplicate the data set twice for each mediator,

mydata1<-mydata
mydata2<-mydata
mydata1$Astar<-mydata$A
mydata2$Astar<-1-mydata$A

#thus creating 8 (4x2) replicates of original data,
#we name this duplicated data as newmydata;

newmydata<-rbind(mydata1,mydata2,mydata3,mydata4,mydata5,mydata6,mydata7,mydata8)

#data will now look like
> head(newmydata,n=16)
#Data set names
#A: Exposure; Y: Outcome Astar: Duplicated exposure
#Id: Respondent ID; W: Weights untruncated;
#W95: Weights truncated at 95%Percentile.
newmydata[1097:1112,]

Data could not be published due to confidentiality reasons

#95th percentile are truncated.
#Exposure weights before truncation
> summary(newmydata[,"W"])
   Min. 1st Qu.  Median    Mean 3rd Qu.    Max.
    0.3     1.0     1.1     2.0     1.5    53.6
#Exposure weights after truncation
> summary(newmydata[,"W95"])
   Min. 1st Qu.  Median    Mean 3rd Qu.    Max.
   0.31    1.03    1.12    1.75    1.46   22.17
\end{verbatim}

\subsection*{R-code for binary exposure, multinomial mediator and binary outcome}
This example is similar to that described in main document. Here we use socio-economic position as the exposure (denoted $A$), measured at 9 months of (child) age, the mediator maternal psychological distress (denoted by $M_1, M_2, M_3,\text{and} M_4$), measured at ages 3, 5, 7 and 11, and the outcome, child development (denoted $Y$), measured at age 11. The considered confounders of the exposure-outcome relation are maternal education, ethnicity and marital status as a proxy for family environment. The time-varying confounders include the number of siblings, the previous outcomes and mediators measured at ages 3, 5 and 7.
\\
\textbf{Step-1:} We first load the data set (e.g. mental health data) into a vector named ``mydata'' in R. Then fit a generalised linear regression to the mediator (\texttt{mental psychological distress}) at each of the four time points separately, conditioning on the exposure \texttt{A} and the confounders listed above (which include the history of outcomes and mediators at each time). Since the mediator is a categorical variable we used the vector generalised linear model (\texttt{vglm}), with family being multinomial, from the VGAM library in R.
\small{
\begin{verbatim}
m1creg <- vglm(M1~A+factor(L01)+factor(L02)+factor(L03)+factor(L04)+factor(L05)
                    ,data=mydata,family=multinomial)
m2creg <- vglm(M2~factor(M1)+A+Y1+factor(L01)+factor(L02)+factor(L03)
               +factor(L04)+factor(L05),data=mydata,family=multinomial)
m3creg <- vglm(M3~factor(M2)+factor(M1)+Y2+Y1+A+factor(L01)+factor(L02)
               +factor(L03)+factor(L04)+factor(L05),data=mydata,family=multinomial)
m4creg <- vglm(M4~factor(M3)+factor(M2)+factor(M1)+Y2+Y1+Y3+A+factor(L01)
               +factor(L02)+factor(L03)+factor(L04)+factor(L05),
               data=mydata,family=multinomial)
\end{verbatim}
}

\textbf{Step-2:} In this step we create a new data set by replicating the original two times. Next we create a new variable \texttt{Id}, which is a subject identifier, and \texttt{Astar} which takes once takes the value contained in \texttt{A} for each subject, and once the opposite value \texttt{1-A}.
\begin{verbatim}
N <- nrow(mydata)
mydata$id <- 1:N
mydata1 <- mydata
mydata2 <- mydata
mydata1$Astar <- mydata$A
mydata2$Astar <- 1-mydata$A
newmydata <- rbind(mydata1,mydata2)
\end{verbatim}

\textbf{Step-3:} The weights are computed from the predicted probabilities of the above multinomial regressions. For the numerator computation in weight equation 11  we used $1-A$ and for the denominator we use used $A$ as explanatory variables:
\begin{verbatim}
#denominator weights
m1pdr <- as.matrix(predict(m1reg,type="response",newdata=mydata))
                    [cbind(1:nrow(mydata),mydata$M1)]
m2pdr <- as.matrix(predict(m2reg,type="response",newdata=mydata))
                   [cbind(1:nrow(mydata),mydata$M2)]
m3pdr <- as.matrix(predict(m3reg,type="response",newdata=mydata))
                    [cbind(1:nrow(mydata),mydata$M3)]
m4pdr <- as.matrix(predict(m4reg,type="response",newdata=mydata))
                   [cbind(1:nrow(mydata),mydata$M4)]

#numerator weights
Dattemp <- mydata
Dattemp$A <- 1 - Dattemp$A
m1pnr<-as.matrix(predict(m1reg,type="response",newdata=Dattemp))
                   [cbind(1:nrow(mydata),mydata$M1)]
m2pnr<-as.matrix(predict(m2reg,type="response",newdata=Dattemp))
                   [cbind(1:nrow(mydata),mydata$M2)]
m3pnr<-as.matrix(predict(m3reg,type="response",newdata=Dattemp))
                  [cbind(1:nrow(mydata),mydata$M3)]
m4pnr<-as.matrix(predict(m4reg,type="response",newdata=Dattemp))
                  [cbind(1:nrow(mydata),mydata$M4)]
\end{verbatim}
The predict function produces a matrix with probabilities for all 9 possible values of the mediator. By including the argument$[cbind(1:nrow(mydata),mydata\$M1)]$ we select the probabilities corresponding to the values actually observed for the mediator. This trick requires that the levels of the mediators start from 1 in the data set. Since in our data set the mediator values start from 0 we changed the original values of mediators  $mydata\$M <- mydata\$M+1$, so that the above command($cbind(1:nrow(mydata)$) works.

\textbf{Step-4:} The weights for the exposure were created using the logistic regression conditioned on the baseline confounders.
\begin{verbatim}
#Unstabilized weights
Aregdr <- glm(A~factor(L01)+factor(L02)+factor(L03)+factor(L04)+factor(L05)
                    ,data=mydata,family=binomial)
ap <- predict(Areg,type="response")
wa <- ifelse(mydata$A==1,1/ap,1/(1-ap)) #weight of A.
\end{verbatim}
\textbf{Step-5:} Finally the logistic natural effects model for the dichotomous outcome can be fitted. To obtain sandwich standard errors that account for the repeated measures nature of the data, we use the \texttt{geeglm} function from the package \texttt{geepack} using the \texttt{Id} variable to indicate dependence. Note that the \texttt{geeglm} function requires that the observations be sorted by the \texttt{Id} variable.
\begin{verbatim}
library(geepack}
newmydata <- newmydata[order(newmydata$id),]
fit_out <- geeglm(Y~A+Astar+t+I(t*A)+I(T*Astar),data=newmydata,corstr="independence",
                       family="binomial", weights=W,id=newmydata$id,scale.fix=T)
> summary(fit_out)

Call:
geeglm(formula = Y ~ A + Astar + T + I(T * A) + I(T * Astar),
    family = "binomial", data = gdat, weights = W, id = newmydata$id,
    corstr = "independence", scale.fix = T)

 Coefficients:
             Estimate  Std.err   Wald Pr(>|W|)
(Intercept)  -2.10513  0.11583 330.28  < 2e-16 ***
A             0.24105  0.16272   2.19  0.13851
Astar         0.11682  0.01725  45.87  1.3e-11 ***
T            -0.13763  0.03981  11.95  0.00055 ***
I(T * A)      0.14047  0.06084   5.33  0.02097 *
I(T * Astar)  0.00609  0.01133   0.29  0.59085
---
Signif. codes:  0 ‘***’ 0.001 ‘**’ 0.01 ‘*’ 0.05 ‘.’ 0.1 ‘ ’ 1

Scale is fixed.

Correlation: Structure = independence
Number of clusters:   5189   Maximum cluster size: 8
\end{verbatim}
\subsubsection*{Code for computing the variance}
The below code is for computing the \textbf{variance} using the method described in section-4 of the manuscript.\\

\subsubsection*{Code for drawing samples}
This section has three components; 1) adding error to coefficients of exposure regression, 2) adding error to the coefficients of mediator regression and 3) computing probabilities of the multinomial regression
\begin{verbatim}
#Part-1: generating new set of weights for the exposure regression
library(mvtnorm)
awts<-function(regfit){
sreg<-summary(regfit)
cv<-sreg$cov.unscaled
mdl<-model.matrix(regfit)
cf<-sreg$coefficients[,1]
l<-length(cf)
pe<-rmvnorm(1,rep(0,l),cv)
ncf<-cf+pe
ncf<-t(ncf)
nodds<-mdl%*%ncf
nweights<-exp(nodds)/(1+exp(nodds))
return(nweights)
}

#Part-2: generating the error to be added to coefficients for the mediator regression
library(mvtnorm)
mcoeffs<-function(mfit){
coeff<-coef(mfit,matrix=TRUE)
cv<-vcov(mfit)
l<-dim(coeff)[1]
k<-dim(coeff)[2]
m<-l*k
#Since the mfit has coefficients corresponding to every level of a multinomial mediator
#information from the coeff and cv matrix above need to be extracted carefully
#it is for this reason we at first developed an index, s,
#this index s now allows to extract the correct set of coefficients  and covariances
#corresponding to a level of a mediator.
s<-seq(1,m,k)
#s: 1  9 17 25 33 41 49 57 65 73 81 89 97
ncoeff<-matrix(0,k,l)
cf<-coeff[,1]
cv1<-as.matrix(cv[s,s])
me<-rmvnorm(1,rep(0,l),cv[s,s])
ncoeff[1,]<-coeff[,1]+me
for(i in 1:(k-1)){
cf<-coeff[,i+1]
cv1<-as.matrix(cv[s+i,s+i])
men<-rmvnorm(1,rep(0,l),cv1)
ncoeff[i+1,]<-cf+men
}
ncoeff<-t(ncoeff)
return(ncoeff)
}
#Part-3 program for computing the predicted probabilities in case of multinomial mediator
trial<-function(mdl,ncoef){
ro<-dim(mdl)[1]
co<-dim(mdl)[2]
k<-ncol(ncoef)
n<-ro/k
rsm<-matrix(0,n,k)
fs<-seq(1,ro,k)
ss<-seq(1,co,k)
rsm[,1]<-mdl[fs,ss]%*%ncoef[,1]
for(i in 1:(k-1)){
pv<-mdl[fs+i,ss+i]%*%ncoef[,(i+1)]
rsm[,(i+1)]<-pv
}
rsm<-exp(rsm)
rs<-rowSums(rsm)
s<-matrix(0,n,k)
for(i in 1:n){
s[i,]<-rsm[i,]/(1+rs[i])
}
rss<-rowSums(s)
s<-as.data.frame(s)
s$V9<-1-rss
s<-as.matrix(s)
return(s)
}
\end{verbatim}
\subsubsection*{Code for computing the variance}
Once the weights are recomputed from each simulation using the above code then the outcome regression needs to be performed for each simulation. Outcome regression coefficients and the \textbf{variance}-covariance matrix from each of this simulation are stored to compute the final variance estimate that accounts for varying weights. To do the variance computation using the stored estimates we used the following code.
\begin{verbatim}
#Code for exposure regression before conducting simulations
Aregdr<-glm(Atemp~factor(L01)+factor(L02)+factor(L03)
            +factor(L04)+factor(L05),data=mydata,family=binomial)
#code for mediator regression from each time point before simulations
m1reg<-vglm(M1~Atemp+factor(L01)+factor(L02)
            +factor(L03)+factor(L04)+factor(L05)
            ,data=mydata,family=multinomial)
m2reg<-vglm(M2~factor(M1)+Atemp+Y1+factor(L01)
            +factor(L02)+factor(L03)+factor(L04)+factor(L05)
            ,data=mydata,family=multinomial(),maxit=500)
m3reg<-vglm(M3~factor(M2)+factor(M1)+Y2+Y1+Atemp+factor(L01)
            +factor(L02)+factor(L03)+factor(L04)
            +factor(L05),data=mydata,family=multinomial())
m4reg<-vglm(M4~factor(M3)+factor(M2)+factor(M1)+Y2+Y1+Y3+Atemp
            +factor(L01)+factor(L02)+factor(L03)
            +factor(L04)+factor(L05),data=mydata,family=multinomial())
mdl1<-model.matrix(m1reg)
mdl2<-model.matrix(m2reg)
mdl3<-model.matrix(m3reg)
mdl4<-model.matrix(m4reg)

mydata$Atemp<-1-mydata$A
m1creg<-vglm(M1~Atemp+factor(C1)+factor(L01)+factor(L02)+factor(L03)
             +factor(L04)+factor(L05),data=mydata,family=multinomial)
m2creg<-vglm(M2~factor(M1)+Atemp+Y1+factor(L01)+factor(L02)+factor(L03)
             +factor(L04)+factor(L05),data=mydata,family=multinomial(),maxit=500)
m3creg<-vglm(M3~factor(M2)+factor(M1)+Y2+Y1+Atemp+factor(L01)+factor(L02)+factor(L03)
             +factor(L04)+factor(L05),data=mydata,family=multinomial())
m4creg<-vglm(M4~factor(M3)+factor(M2)+factor(M1)+Y2+Y1+Y3+Atemp+factor(L01)
             +factor(L02)+factor(L03)+factor(L04)+factor(L05)
             ,data=mydata,family=multinomial())
mdlc1<-model.matrix(m1creg)
mdlc2<-model.matrix(m2creg)
mdlc3<-model.matrix(m3creg)
mdlc4<-model.matrix(m4creg)

M<-500
t<-matrix(c(0,1,2,3),1,4)
#full code of simulations using above code
for(i in 1:M){
apdr<-awts(Aregdr)
wa<-ifelse(mydata$A==1,1/apdr,(1-mydata$A)/(1-apdr)) #weight of A.
wa_95<-ifelse(wa>=quantile(wa,.95),quantile(wa,.95),wa)


c1<-mcoeffs(m1reg)
c2<-mcoeffs(m2reg)
c3<-mcoeffs(m3reg)
c4<-mcoeffs(m4reg)
mt1nr<-trial(mdlc1,c1)[cbind(1:nrow(mydata),mydata$M1)]
mt2nr<-trial(mdlc2,c2)[cbind(1:nrow(mydata),mydata$M2)]
mt3nr<-trial(mdlc3,c3)[cbind(1:nrow(mydata),mydata$M3)]
mt4nr<-trial(mdlc4,c4)[cbind(1:nrow(mydata),mydata$M4)]

mt1dr<-trial(mdl1,c1)[cbind(1:nrow(mydata),mydata$M1)]
mt2dr<-trial(mdl2,c2)[cbind(1:nrow(mydata),mydata$M2)]
mt3dr<-trial(mdl3,c3)[cbind(1:nrow(mydata),mydata$M3)]
mt4dr<-trial(mdl4,c4)[cbind(1:nrow(mydata),mydata$M4)]

#creation of new weights
wM1<-mt1nr/mt1dr
mydata$wM1<-wM1*wa

wM2<-(mt2nr/mt2dr)*wM1
mydata$wM2<-wM2*wa

wM3<-(mt3nr/mt3dr)*wM2
mydata$wM3<-wM3*wa

wM4<-(mt4nr/mt4dr)*wM3
mydata$wM4<-wM4*wa
#truncated weights
wM1<-mt1nr/mt1dr
mydata$wM1_95<-wM1*wa_95

wM2<-(mt2nr/mt2dr)*wM1
mydata$wM2_95<-wM2*wa_95

wM3<-(mt3nr/mt3dr)*wM2
mydata$wM3_95<-wM3*wa_95

wM4<-(mt4nr/mt4dr)*wM3
mydata$wM4_95<-wM4*wa_95




#creation of final data for analysis
gdat1<-cbind(A=mydata$A,Y=mydata$Y1,Astar=mydata$A,T=0,
             id=1:nrow(mydata),W=wa,W_95=wa_95)
gdat2<-cbind(A=mydata$A,Y=mydata$Y1,Astar=1-mydata$A,T=0,
             id=1:nrow(mydata),W=mydata$wM1,W_95=mydata$wM1_95)
gdat3<-cbind(A=mydata$A,Y=mydata$Y2,Astar=mydata$A,T=1,
             id=1:nrow(mydata),W=wa,W_95=wa_95)
gdat4<-cbind(A=mydata$A,Y=mydata$Y2,Astar=1-mydata$A,T=1,
             id=1:nrow(mydata),W=mydata$wM2,W_95=mydata$wM2_95)
gdat5<-cbind(A=mydata$A,Y=mydata$Y3,Astar=mydata$A,T=2,
             id=1:nrow(mydata),W=wa,W_95=wa_95)
gdat6<-cbind(A=mydata$A,Y=mydata$Y3,Astar=1-mydata$A,T=2,
             id=1:nrow(mydata),W=mydata$wM3,W_95=mydata$wM3_95)
gdat7<-cbind(A=mydata$A,Y=mydata$Y4,Astar=mydata$A,T=3,
             id=1:nrow(mydata),W=wa,W_95=wa_95)
gdat8<-cbind(A=mydata$A,Y=mydata$Y4,Astar=1-mydata$A,T=3,
             id=1:nrow(mydata),W=mydata$wM4,W_95=mydata$wM4_95)

#combining the data
newmydata<-as.data.frame(rbind(gdat1,gdat2,gdat3,gdat4
                            ,gdat5,gdat6,gdat7,gdat8))
newmydata<-newmydata[order(newmydata$id),]

newmydata<-newmydata[order(newmydata$id),]
library(geepack)
fitsepb<-geeglm(Y~A+Astar+T+I(T*A)+I(T*Astar),id=newmydata$id,
                corstr="independence",family="binomial",
                weight=W,data=newmydata,scale.fix=T)

#extracting coefficients and variance covariance from each simulation
cf<-summary(fitsepb)$coefficients
cv<-vcov(fitsepb)
#DE matrix of direct effects
#IDE matrix of indirect effects
$VC variance covariance terms required
DE[i,]<-c(cf[2,1],cf[5,1])
IDE[i,]<-c(cf[3,1],cf[6,1])
VC[i,]<-c(cv[2,2],cv[3,3],cv[5,5],cv[6,6],
          cv[2,3],cv[2,5],cv[2,6],cv[3,5]
          ,cv[3,6],cv[5,6])


#truncated 95
fitsepb_95<-geeglm(Y~A+Astar+T+I(T*A)+I(T*Astar),id=newmydata$id,
                   corstr="independence",family="binomial",
                   weight=W_95,data=newmydata,scale.fix=T)
#extracting coefficients and variance covariance matrices
cf<-summary(fitsepb_95)$coefficients
cv<-vcov(fitsepb_95)

#DE direct effect coefficients matrix corresponding to truncated data
#IDE indirect effect matrix
#Variance covariances corresponding to truncated data model fit.
DE_95[i,]<-c(cf[2,1],cf[5,1])
IDE_95[i,]<-c(cf[3,1],cf[6,1])
VC_95[i,]<-c(cv[2,2],cv[3,3],cv[5,5],cv[6,6]
             ,cv[2,3],cv[2,5],cv[2,6],
              cv[3,5],cv[3,6],cv[5,6])
}

#compute the final estimate of direct, indirect and total effect for each simulation.
#data for these are from the changing weights
#DE:  Stores the coefficients required for computing the direct effect estimates
#IDE: Stores the coefficients required for computing the indirect effect estimates
#TE: total Effects
#DE[,1]: Has the values corresponding to alpha_1 from simulations
#DE[,2]: Has the coefficients corresponding to alpha_4
#IDE[,1]: Has the values corresponding to alpha_2 from simulations
#IDE[,2]: Has the values corresponding to alpha_5
DE_final<-DE[,1]+matrix(DE[,2],M,1)%*%t
IDE_final<-IDE[,1]+matrix(IDE[,2],M,1)%*%t
TE_final<-DE_final+IDE_final

#Step-2
#compute the variances of the DE, IDE and TE for all time points
#data for these are from the variance and covariance extract
#from each simulation with new weights
#VC: Variance covariance estimates stored from each simulation
DEV<-VC[,1]+matrix(VC[,3],M,1)%*%(t^2)+matrix(VC[,6],M,1)%*%(2*t)
IDEV<-VC[,2]+matrix(VC[,4],M,1)%*%(t^2)+matrix(VC[,9],M,1)%*%(2*t)
TEV<-VC[,1]+VC[,2]+matrix(VC[,3],M,1)%*%(t^2)+matrix(VC[,6],M,1)%*%(2*t)
           +matrix(VC[,4],M,1)%*%(t^2)+matrix(VC[,9],M,1)%*%(2*t)
           +matrix(VC[,5],M,1)%*%(2*t)+matrix(VC[,7],M,1)%*%(2*t)
           +matrix(VC[,8],M,1)%*%(2*t)+matrix(VC[,10],M,1)%*%(2*t^2)

#Step-3
#Computing the new variance of the direct, indirect, and total effects
#that account for the weight changes
DEvar<-apply(DE_final,2,var)+apply(DEV,2,mean)
IDEvar<-apply(IDE_final,2,var)+apply(IDEV,2,mean)
TEvar<-apply(TE_final,2,var)+apply(TEV,2,mean)

#Step-4 computing the lower and upper bounds of confidence interval
#which accounts for variance due to changing weights.
#In this step I am using the initial estimates computed from the observed data and
#stored in a matrix labelled "final".
#DEL: Direct effect lower; DEU: Direct effect upper; IDEL: Indirect effect lower;
#IDEU: Indirect effect upper
#TEL: Total effect lower; TEU: Total effect upper.
DEL<-final[1,]-sqrt(DEvar)*1.96
DEU<-final[1,]+sqrt(DEvar)*1.96

IDEL<-final[2,]-sqrt(IDEvar)*1.96
IDEU<-final[2,]+sqrt(IDEvar)*1.96

TEL<-final[3,]-sqrt(TEvar)*1.96
TEU<-final[3,]+sqrt(TEvar)*1.96
#Combining all the estimates to be reported
finalest<-rbind(DEL,final[1,],DEU,IDEL,final[2,],IDEU,TEL,final[3,],TEU)
#the above process is then repeated for the truncated estimates.
\end{verbatim}
If one estimates the total effect using an MSM by setting $A=a^*$ and use the default GEE Std.err, and compare it to the Std.err obtained using above described method they should more or less agree.
\section*{Appendix 3}
\subsection*{Variance estimation}
Let $\hat{\alpha}$ be the estimator of the coefficient $\alpha$ indexing the natural effect model and $\hat{w}$ be the estimated weight.
Then it follows by the law of iterated variance that
\[Var(\hat{\alpha})=E\left\{Var(\hat{\alpha}|\hat{w})\right\}+Var\left\{E(\hat{\alpha}|\hat{w})\right\}.\]
Writing $\hat{\alpha}=\hat{\alpha}(O,\hat{w})$ as a function of the observed data $O$ and the estimated weights $\hat{w}$, we have that
\[E(\hat{\alpha}|\hat{w})=\int \hat{\alpha}(O,\hat{w})f(O|\hat{w})dO=\int \hat{\alpha}(O,\hat{w})\left\{f(O)+o_p(1)\right\}dO,\]
where $o_p(1)$ is a term that converges to zero in probability. That $f(O|\hat{w})=f(O)+o_p(1)$ can be seen upon noting that $\hat{w}$ converges to a deterministic function of $O$ in probability.
It follows that $E(\hat{\alpha}|\hat{w})$ equals the mean of $\hat{\alpha}(O,w)$ with $w$ substituted by $\hat{w}$, up to an $o_p(1)$ term.
Likewise, $Var(\hat{\alpha}|\hat{w})$ equals the variance of $\hat{\alpha}(O,w)$ with $w$ substituted by $\hat{w}$, up to an $o_p(1)$ term.
The variance of $\hat{\alpha}(O,w)$ can be consistently estimated using a sandwich estimator, considering the weight $w$ as known.
It follows that $E\left\{Var(\hat{\alpha}|\hat{w})\right\}$ can be consistently estimated as
\[\frac{1}{B}\sum_{j=1}^B Var(\hat{\alpha}|\hat{w}_j).\]
Further, the mean $E(\hat{\alpha}|\hat{w})$ can be consistently estimated as $\hat{\alpha}(O,\hat{w})$. It follows that $Var\left\{E(\hat{\alpha}|\hat{w})\right\}$ can be consistently estimated as
\[\frac{1}{B-1}\sum_{j=1}^B \left(\hat{\alpha}^{(j)} - B^{-1}\sum_{k=1}^B \hat{\alpha}^{(k)}\right)\left(\hat{\alpha}^{(j)} - B^{-1}\sum_{k=1}^B \hat{\alpha}^{(k)}\right)'.\]
\subsection*{Simulations for standard error estimation}
In this section we present the code for simulating the data that was used for conducting single mediation analysis using natural effect models. We also present the code that was used in estimating the standard errors. Simulations were performed on a sample of 1000 observations and with 1000 perturbed bootstrap estimates.
\begin{verbatim}
#program for checking the variance.
#I am generating the values of exposure, mediator, confounders and outcome as Binary
n<-1000
set.seed(562)
C1<-rnorm(n)
C2<-rbinom(n,1,0.4)
A<-rbinom(n,1,(exp(0.05+0.1*C1+0.2*C2)/(1+exp(0.05+0.1*C1+0.2*C2))))
M<-rbinom(n,1,(exp(0.05+0.1*A-0.1*C1-0.2*C2)/(1+exp(0.05+0.1*A-0.1*C1-0.2*C2))))
Y<-rbinom(n,1,(exp(0.05+0.1*A+0.1*M-0.1*C1-0.2*C2)/(1+exp(0.05+0.1*A+0.1*M-0.1*C1-0.2*C2))))


#storing of the data set
sdat<-as.data.frame(cbind(Y,M,A,C1,C2))
write.table(sdat,file="sdat.csv",sep=",",row.names=FALSE)


#Conducting natural effects model in single mediator case.
Atemp<-A
sdat$Atemp<-A
#creation of Inverse probability treatement weights
Aregdr<-glm(Atemp~C1+C2,data=sdat,family=binomial)
ap<-fitted.values(Aregdr,type="response")
wa<-ifelse(sdat$A==1,1/ap,1/(1-ap)) #weight of A.
#Creation of weights of M for eachi time period
M1reg<-glm(M~Atemp+C1+C2,data=sdat,family=binomial("logit"))
sdat$Atemp<-sdat$A
tdr1<-predict.glm(M1reg,newdata=sdat,type="response")
sdat$Atemp<-1-sdat$A
tnr1<-predict.glm(M1reg,newdata=sdat,type="response")

#creation of mediator weights
W<-tnr1/tdr1

#Adding Weight of M to data set
sdat$wM1<-(tnr1/tdr1)*wa

#expanding data
gdat1<-cbind(A=sdat$A,Y=sdat$Y,Astar=sdat$A,id=1:nrow(sdat),W=wa)
gdat2<-cbind(A=sdat$A,Y=sdat$Y,Astar=1-sdat$A,id=1:nrow(sdat),W=sdat$wM1)


#combining the data
gdat<-as.data.frame(rbind(gdat1,gdat2))
gdat<-gdat[order(gdat$id),]

#conducting final estimation using the GEEGLM model.
library(geepack)
library(BSagri)
fitsepb<-geeglm(Y~A+Astar,id=gdat$id,corstr="independence",
                family="binomial",weight=W,data=gdat,scale.fix=T)
exp(fitsepb$coefficients)
de_model<-fitsepb$coefficients[,2]
ide_model<-fitsepb$coefficients[,3]
de_odds<-exp(de_model)
ide_odds<-exp(ide_model)
te<-de_model+ide_model
te-odds<-exp(te)
vcc<-vcov(fitsepb)
te_v<-vcc[2,2]+vcc[3,3]+2*vcc[2,3]
te_se<-sqrt(te_v)


#Perturb method of bootstrap.
#program for the weight generation.
#generating new set of weights for the exposure without re-sampling the data
awts<-function(regfit){
sreg<-summary(regfit)
cv<-sreg$cov.unscaled
mdl<-model.matrix(regfit)
library(mvtnorm)
cf<-sreg$coefficients[,1]
l<-length(cf)
pe<-rmvnorm(1,rep(0,l),cv)
ncf<-cf+pe
ncf<-t(ncf)
nodds<-mdl%*%ncf
nweights<-exp(nodds)/(1+exp(nodds))
return(nweights)
}

#genrating the distributions of coefficients for the mediator data
#generating new set of weights for the exposure without re-sampling the data
mwts<-function(regfit,Atemp){
sreg<-summary(regfit)
cv<-sreg$cov.unscaled
mdl<-as.data.frame(model.matrix(regfit))
mdl$Atemp<-Atemp
mdl<-as.matrix(mdl)
library(mvtnorm)
cf<-sreg$coefficients[,1]
l<-length(cf)
pe<-rmvnorm(1,rep(0,l),cv)
ncf<-cf+pe
ncf<-t(ncf)
nodds<-mdl%*%ncf
nweights<-exp(nodds)/(1+exp(nodds))
return(nweights)
}


#using 1000 perturb samples
M<-1000
Dat<-sdat
DE<-matrix(0,M,1)
IDE<-matrix(0,M,1)
VC<-matrix(0,M,3)
library(geepack)
library(BSAgri)
#Final estimates
for(i in 1:M){
apdr<-awts(Aregdr)
wa<-ifelse(Dat$A==1,1/apdr,1/(1-apdr)) #weight of A.
Atemp<-Dat$A
wdr<-mwts(M1reg,Atemp)
Atemp<-1-Dat$A
wnr<-mwts(M1reg,Atemp)
wM1<-wnr/wdr
Dat$wM1<-wM1*wa

#creation of final data for analysis
gdat1<-cbind(A=Dat$A,Y=Dat$Y,Astar=Dat$A,id=1:nrow(Dat),W=wa)
gdat2<-cbind(A=Dat$A,Y=Dat$Y,Astar=1-Dat$A,id=1:nrow(Dat),W=Dat$wM1)

#combining the data
gdat<-as.data.frame(rbind(gdat1,gdat2))
gdat<-gdat[order(gdat$id),]


fitsepb<-geeglm(Y~A+Astar,id=gdat$id,corstr="independence",
                family="binomial",weight=W,data=gdat,scale.fix=T)
cf<-summary(fitsepb)$coefficients
cv<-vcov(fitsepb)

DE[i,]<-c(cf[2,1])
IDE[i,]<-c(cf[3,1])
VC[i,]<-c(cv[2,2],cv[3,3],cv[2,3])
}

#computing the variance of DE using the perturb bootstrap
va<-var(DE)
mvde<-mean(VC[,1])
tvde<-va+mvde

#computing the variance of IDE using the perturb bootstrap
vide<-var(IDE)
mvide<-mean(VC[,2])
tvide<-vide+mvide

exp(mean(DE))
exp(mean(IDE))

#computing the total effect and its variance.
Te<-DE+IDE
exp(mean(Te))
vte<-VC[,1]+VC[,2]+2*VC[,3]

#end of program.
\end{verbatim}
\begin{table}
\centering{
\caption{Odds ratio (OR) estimates of direct and indirect effects estimated using simulated data. These estimates are computed using natural effect models. The standard error (Std. err) estimates are obtained using empirical estimates, generalized estimating equations (GEE) and perturbed bootstrap.}\label{wt4}
\begin{center}
\begin{tabular}{|p{1.5cm}|p{1.2cm}|p{1.2cm}|p{1.2cm}|p{1.2cm}|p{1.2cm}|p{1.2cm}|}
 \hline
 Estimate&\multicolumn{2}{c|}{Empirical}&\multicolumn{2}{c|}{GEE} &\multicolumn{2}{c|}{Perturbed Bootstrap}\\
 \hline
         &OR&Std. err&OR&Std. err&OR&Std. err\\
         \hline
 Direct &1.01&0.124&1.01&0.127&1.01&0.127\\
 Indirect &0.99&0.004&0.99&0.001&0.99&0.008\\
  \hline
\end{tabular}
\end{center}}
\end{table}
\end{appendices}
\end{document}